\documentclass{article}

\usepackage{amsmath}
\usepackage{amssymb}
\usepackage{url}
\usepackage{graphicx}

\newtheorem{theorem}{Theorem}

\title{Boundary slopes of some non-Montesinos knots}
\author{Joshua Howie\\University of Melbourne}

\begin{document}

\maketitle

It is shown that there exist alternating non-Montesinos knots whose essential spanning surfaces with maximal and minimal boundary slopes are not realised by the checkerboard surfaces coming from a reduced alternating planar diagram.

\subsection*{Slope Diameter}

The boundary slope diameter of a knot $K$ is
\[d(K) = \max\{|s-s'|:s,s'\in\mathcal{B}(K)\setminus\{\infty\}\},\]
where $\mathcal{B}(K)$ is the finite set of boundary slopes for $K$.
Let $c(K)$ denote the crossing number of $K$.

\begin{theorem}[Ichihara-Mizushima~\cite{ichmiz1}]
If $K$ is a Montesinos knot, then
\[ d(K)\leq2c(K),\]
with equality if $K$ is alternating and Montesinos.
\end{theorem}

This theorem was proved for the case of two-bridge knots in~\cite{mmr2br}.

A spanning surface for a knot $K$ is a surface $\Sigma$ embedded in $S^3$ such that $\partial\Sigma = K$. Let $X$ denote the knot exterior, this is the closure of $S^3\setminus N(K)$. The slope of a spanning surface $\Sigma$ is the slope of the surface $\Sigma\cap X$ on $\partial X$. 

Let $\mathcal{S}(K)$ be the subset of $\mathcal{B}(K)$ that consists of boundary slopes that are realised by spanning surfaces. Let $d_{\mathcal{S}}(K)$ be the diameter of $\mathcal{S}(K)$. Note that $\mathcal{S}(K)\subset 2\mathbb{Z}$, whereas $\mathcal{B}(K)\subset\mathbb{Q}\cup\{\infty\}$.

\begin{theorem}[Curtis-Taylor~\cite{curttay}]
If $K$ is an alternating knot, then
\[ d_{\mathcal{S}}(K) = 2c(K).\]
\end{theorem}

It is the purpose of this note to provide a counterexample to Theorem 2.

The proof given in~\cite[page 1350]{curttay} contains a mistake, where a result of Adams and Kindred~\cite{adkind} has been incorrectly stated. Adams and Kindred showed that given an essential spanning surface $\Sigma$ for an alternating knot $K$, then a spanning surface $\Sigma'$ can be obtained from a basic layered surface $S$ for $K$, by adding some number of handles or crosscaps to $S$, such that $\Sigma'$ has the same orientability, slope and genus as $\Sigma$. 

Curtis and Taylor have quoted this result as: given an essential spanning surface $\Sigma$ for an alternating knot $K$, then there exists a basic layered surface $S$ which has the same orientability, slope and genus as $\Sigma$. This is not true because adding a crosscap changes the slope by $\pm 2$, so it is possible to obtain spanning surfaces for a knot $K$ whose slopes lie outside the range of slopes of basic layered surfaces. In particular, there may exist essential spanning surfaces for an alternating knot which has a slope bigger or smaller than both the checkerboard surfaces associated to a reduced, alternating diagram. We will show that such surfaces can exist for non-Montesinos knots.

\subsection*{Generalised Alternating Knots}

As introduced by Hayashi~\cite{hayalt} and Ozawa~\cite{ozntriv}, a knot $K$ is generalised alternating if it has a projection onto a closed orientable embedded surface $F$, 
\[\pi:F\times I\rightarrow F,\] 
where $K\subset F\times I\subset S^3$ such that:

\begin{enumerate}

\item $\pi(K)$ is alternating on $F$, and
\item $\pi(K)$ is prime.

\end{enumerate}

A knot projection on to a surface $F$ is prime if given any loop $l\subset F$, such that $|l\cap\pi(K)| = 2$, then $l$ bounds a disk $D\subset F$, such that $D$ contains only a single embedded arc of $\pi(K)$. This definition in equivalent to the usual definition of prime for a diagram on $S^2$, however in the case of a higher genus projection surface $F$, there exist loops $l\subset F$ that do not bound disks.

Generalised alternating implies that the regions of $F\setminus\pi(K)$ are disks and that any essential loop on $F$ meets $\pi(K)$ at least four times. $K$ bounds checkerboard surfaces relative to its projection on $F$, and at least one of these surfaces is non-orientable.

\begin{theorem}[Ozawa~\cite{ozntriv}]
Let $\pi(K)$ be a generalised alternating projection of a knot $K$ onto a closed orientable surface $F$. Then both the checkerboard surfaces relative to $F$ are essential in $X$.
\end{theorem}

The essentiality of checkerboard surfaces is proved for a larger class of alternating surface projections in~\cite{howrub}, which includes all generalised alternating projections.

A method for enumerating generalised alternating projections onto the torus was developed in~\cite{howie1}. That paper also shows how to construct many examples of generalised alternating projections onto higher genus surfaces.

\subsection*{Counterexamples}

The following examples of generalised alternating projections were found using the methods of~\cite{howie1}. The easiest way to check that they are generalised alternating is to check that all regions of $F\setminus\pi(K)$ are disks, and for each black region, that all of its adjacent white regions are distinct.

\begin{figure}[h]
    \centering
    \includegraphics[scale=0.6, trim = 80mm  12mm 2mm 16mm, clip]{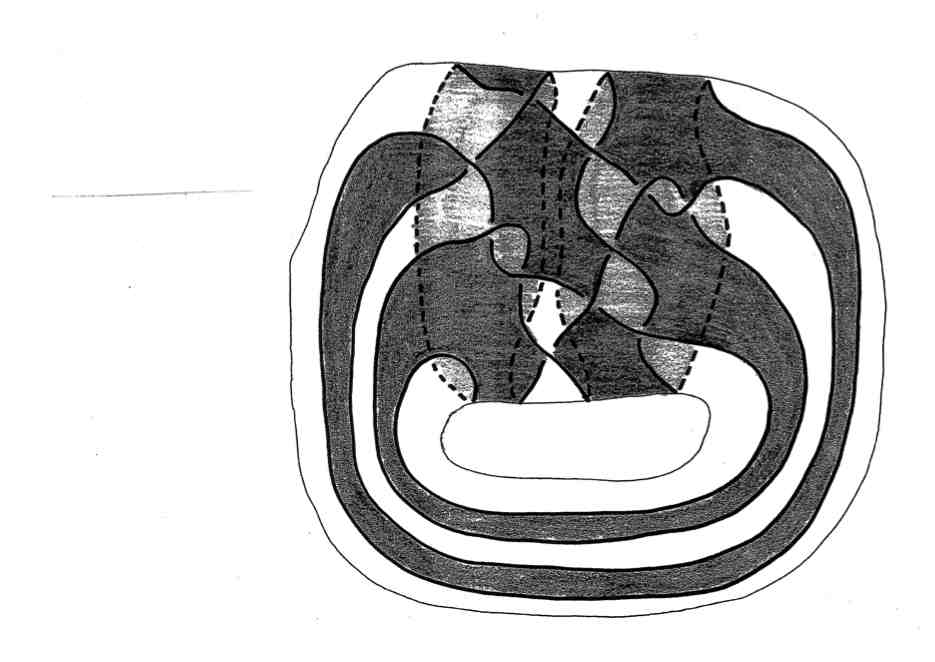}
    \caption{Generalised alternating projection of $8_{17}$.}
    \label{figk1}
\end{figure}

Let $K=8_{17}=8a_{14}$. In~\cite{garoujs}, it was stated that by work of Kabaya~\cite{kabidl},
\[\mathcal{B}(K) = \{-14,-8,-6,-4,-2,0,2,4,6,8,14,\infty\},\]
so clearly $d(K) > 2c(K)=16$. We will show that in fact $d_{\mathcal{S}}(K)= 28 > 2c(K)$.

Figure~\ref{figk1} shows a generalised alternating projection of $K$. In it the black surface has slope $+14$ and the white surface has slope $-8$. A method for calculating the slopes of spanning surfaces is outlined in~\cite{adkind}. The black surface is non-orientable and has Euler charcteristic $-8$.

There is a unique reduced alternating diagram of $8_{17}$ on $S^2$. The two checkerboard surfaces have slopes $+8$ and $-8$. Each surface has Euler characteristic $-3$ and is non-orientable. If we take the planar checkerboard surface with slope $+8$ as our basic layered surface, add three crosscaps in the appropriate way, and add one handle, then we obtain a spanning surface which is non-orientable, has slope $+14$ and Euler charcteristic $-8$, demonstrating that this example does not contradict~\cite{adkind}.

If we reflect the given toroidal projection of $8_{17}$ in the plane, we obtain a different projection of $8_{17}$. This projection is also generalised alternating on a torus. It has checkerboard surfaces with slopes $-14$ and $+8$. Note that $8_{17}$ is amphichiral. Therefore $d_{\mathcal{S}}(K)= 28$ since both the slopes $-14$ and $+14$ are realised by essential spanning surfaces.

In~\cite{howie1}, it is shown that if a knot has diameter $28$, then any generalised alternating projection must be onto a surface of genus at most $2$. It is then natural to ask if there is a generalised alternating projection of $8_{17}$ onto the torus or double torus that has checkerboard surfaces with slopes $+14$ and $-14$?

\begin{figure}[h]
    \centering
    \includegraphics[scale=0.5, trim = 50mm  12mm 2mm 32mm, clip]{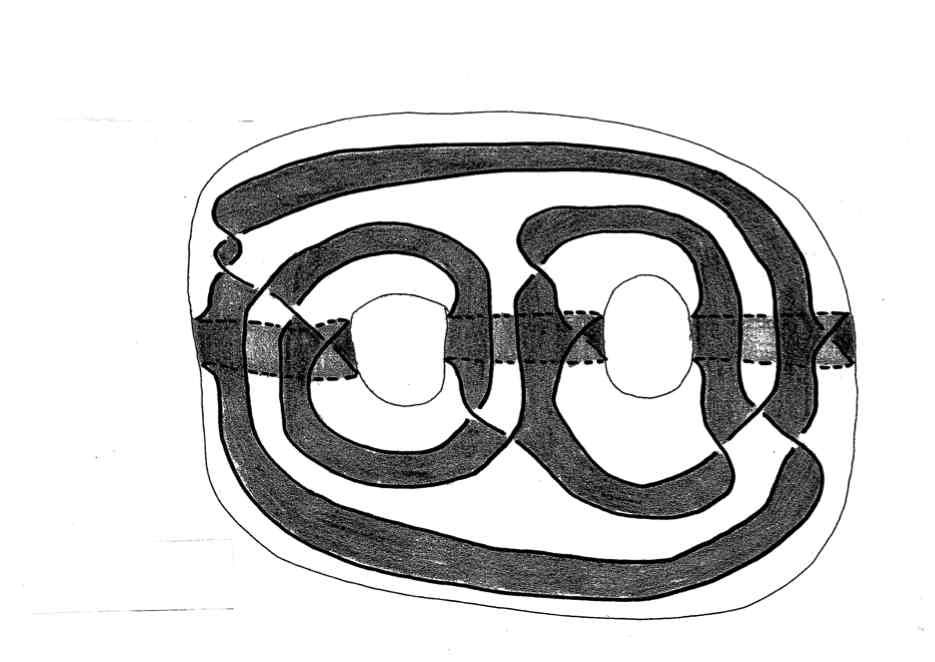}
    \caption{Generalised alternating projection of $12a_{603}$.}
    \label{figk2}
\end{figure}

Figure~\ref{figk2} shows a projection of the knot $12a_{603}$ which is generalised alternating with 13 crossings on the double torus. Therefore the checkerboard surfaces are distance $26$ apart, yet $2c(K)=24$. The black surface has slope $-10$ and the white surface has slope $+16$. This gives a second counterexample to Theorem 2. The slopes of the planar checkerboard surfaces for $12a_{603}$ are $-10$ and $+14$. Adding one crosscap and two handles to the planar checkerboard surface with slope $+14$ produces a spanning surface with the same slope, Euler characteristic and orientability as the white surface pictured.

\begin{figure}[h]
    \centering
    \includegraphics[scale=0.45, trim = 20mm  2mm 2mm 12mm, clip]{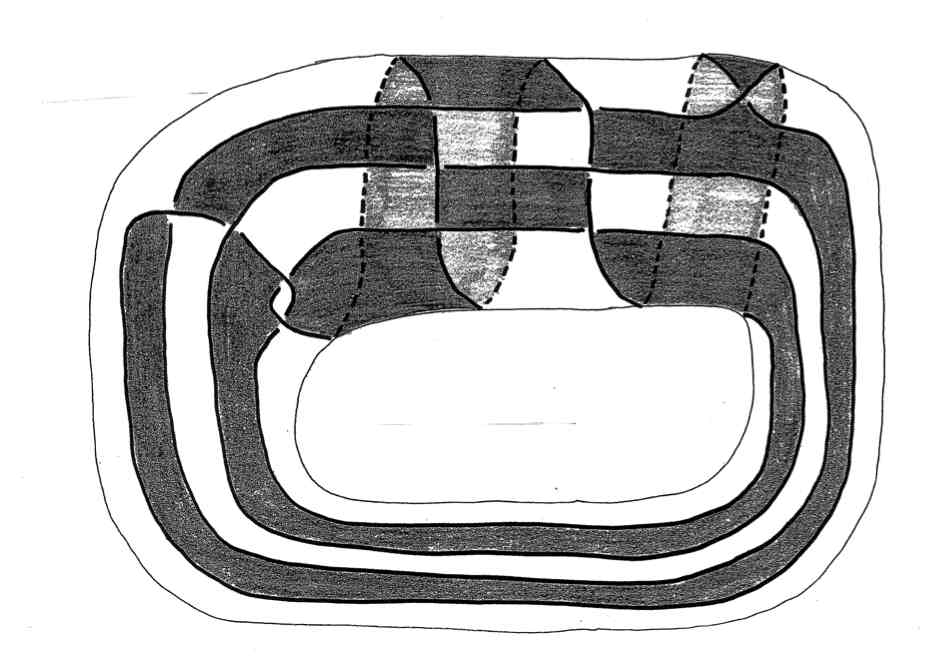}
    \caption{Generalised alternating projection of $10_{161}$.}
    \label{figk3}
\end{figure}

We give one final example, this time of a non-alternating, non-Montesinos knot. Figure~\ref{figk3} shows the knot $K=10_{161}=10n_{31}$. We will show that $d_{\mathcal{S}}(K)\geq 22 > 2c(K)$.
Every knot $K$ has a Newton polygon $N_K$ coming from the $A$-polynomial. This gives a list of visible boundary slopes for $K$. For $K=10_{161}$, see Culler~\cite{cullap}, $N_K$ has $16$ sides, which give five slopes, 
\[\mathcal{B}(K)\supseteq \{-20, -18, -9, -8, 2\}.\]
Of course $0$ is also a slope for $K$.

The black surface has slope $+2$ and the white surface has slope $-20$. Therefore the diameter of $\mathcal{S}(K)$ is at least $22$, whereas its crossing number is only $10$.

Together these three examples establish:

\begin{theorem}
There exist non-Montesinos knots, both alternating and non-alternating, such that
\[ d_{\mathcal{S}}(K) > 2c(K).\]
\end{theorem}

\newpage

\bibliography{largediam}
\bibliographystyle{plain}

\noindent \textsc{Email}: j.howie@student.unimelb.edu.au

\end{document}